# Consensus of fractional-order multi-agent systems: Matrix inequality approach

Elyar Zavary
*Advanced Control Systems Laboratory, School of Electrical and Computer Engineering*
*Tarbiat Modares University*
Tehran, Iran
elyar.zavary@yahoo.com

Mahdi Sojoodi*
*Advanced Control Systems Laboratory, School of Electrical and Computer Engineering*
*Tarbiat Modares University*
Tehran, Iran
sojoodi@modares.ac.ir

*Abstract*—This paper investigates a decentralized dynamic output feedback controller for the robust consensus of fractional-order uncertain multi-agent systems. The procedure is decentralized meaning that each agent can only receive output information from its neighbor agents. By an appropriate state transformation, we reduce the consensus problem to stability one, which leads to a matrix inequality. Additionally, a Homotopy-based method is utilized to solve the resulted matrix inequality and to find the controller's unknown parameters that guarantee the robust consensus of agents. Finally, an illustrative example on multi-agent permanent magnet synchronous motor system and a numerical example are presented to show the effectiveness of the proposed control approach.

*Keywords— Dynamic output feedback, Fractional-order multi-agent systems, Uncertainty, Consensus, Fixed order controller.*

## I. INTRODUCTION

Modeling and study of multi-agent systems have attracted tremendous attention in recent years. This is partly due to their potential applications in many areas, including control theory, mathematics, biology, physics, computer science, and robotics. Consensus is the concept of reaching an agreement considering the states of all agents [1] and plays an important role in multi-agent systems. Examples include flocking [2], formation control [3], cooperative control [4], distributed sensor networks [5], synchronization between the motors, and so on.

The important results of the above literature focus on the consensus problems of multi-agent systems with integer-order dynamical equation such as consensus algorithms of single integrator dynamic systems [7] ,[6] as well as double-integrator dynamics [8],[9] or even high order dynamic systems [10]. Whereas a great number of natural, biological, and practical engineering systems are preferred to be modeled by fractional-order dynamical equations, due to their ability of describing memory and heredity specifications of such systems and phenomena. Synchronized motion of agents in fractional circumstances such as motion of underwater vehicles in lentic lakes and unmanned aerial vehicles in windy and rainy conditions [11], chemotaxis behavior and collective food seeking of microorganisms [12] are examples of these kind of systems and processes. Indeed, integer-order differential equations are incapable of describing dynamics of such systems due to memory and hereditary properties of them. Accordingly, it is meaningful to study the consensus problems of fractional-order systems [11], [13]–[15].

The problem of robust consensus of fractional-order linear multi-agent systems via static feedback was studied in [13], and [16] investigates the distributed containment control of fractional-order uncertain multi-agent systems. Distributed tracking of heterogeneous nonlinear fractional-order multi-agent systems with an unknown leader is studied in [17] via adaptive pinning control. The multi-consensus problem of fractional-order uncertain multi-agent system is converted into the stability problem of fractional-order systems via proper transformation in [11]. Then, the static output feedback controller is utilized for stabilizing the transformed system. It is worth mentioning that controllers, designed based on dynamic feedback, are always preferable to the static ones because of their more effective control performances, moreover the dynamic controller brings about more degree of freedom in achieving control objectives, in comparison with the static controller [18]. In addition, most of mentioned works use state feedback controller and this kind of controllers require all states. On the other hand, in some cases states are inaccessible because of costly implementation or some physical constraints.

High order controllers obtained by most of the controller design methods have expensive implementation procedure, undesirable reliability, high fragility, and numerous maintenance difficulties. Since controller order reduction techniques may deteriorate the closed-loop efficiency, designing directly a low-, fixed-order controller for a system can be helpful [19].

To the best of our knowledge, there is no result on designing dynamic output feedback controller for the consensus of fractional-order multi-agent systems (FOMAS) in the literature. In this paper, by reducing the consensus problem of the FOMAS to stability problem of the augmented system, a fixed-order approach is utilized to gain the above-mentioned advantages of a low-order controller and those of the dynamic output feedback one, simultaneously. Moreover, sufficient conditions for asymptotical consensus of fractional-order multi-agent system are investigated.

The remainder of this paper is organized as follows. In section 2 some preliminaries about fractional calculus, graph theory for describing multi-agent systems, and the Homotopy method are provided. Section 3 establishes the system modeling. In section 4 robust consensus condition of uncertain fractional-order multi-agent systems are derived.

Numerical simulations are presented in section 5 to illustrate the effectiveness of the proposed algorithm. Finally, section 6 concludes the paper.

## II. PRELIMINARIES

In this section, some basic notations of fractional-order systems calculus, concepts and lemmas together with graph theory are presented.

Notations: In this paper $A \otimes B$ denotes the kronecker product of matrices A and B and the symmetric of matrix $M$ will be shown by $sym(.)$ which is defined by $sym(M) = M^T + M$, and also ↑ is the symbol of pseudo inverse.

### A. Fractional calculus

Fractional calculus is a general major in mathematics but we will note some main concepts such as derivative and integral definitions in this section. There are several definitions for fractional-order derivative, since Caputo definitions initial condition is similar to integer order ones as a physical aspect, Caputo definition is used in this with the following definition

$${}_a^C D_t^\alpha = \frac{1}{\Gamma(\bar{n}-\alpha)} \int_a^t (t-\tau)^{\bar{n}-a-1} \left(\frac{d}{d\tau}\right)^{\bar{n}} f(\tau)d\tau,$$

where $\Gamma(\cdot)$ is Gamma function defined by $\Gamma(\epsilon) = \int_0^\infty e^{-t} t^{\epsilon-1} dt$ and $\bar{n}$ is the smallest integer that is equal or greater than $\alpha$.

### B. Algebraic graph theory

A weighted direct graph is shown by $\mathcal{G} = (\mathcal{V}, \varepsilon, \mathcal{A})$, where $\mathcal{V} = \{v_1, \ldots, v_p\}$ is the vertex set, $\varepsilon \subseteq \mathcal{V} \times \mathcal{V}$ is the edge set, and $\mathcal{A} = (a_{ij})_{p \times p}$ is a nonsymmetrical set. A nonzero $\varepsilon_{ij} \in \varepsilon \subseteq \mathcal{V} \times \mathcal{V}$ indicates that agent $j$ receives information from agent $i$ which leads to a corresponding nonzero $a_{ij} \in \mathcal{A} = (a_{ij})_{p \times p}$ and $a_{ij} = 0$ otherwise. Furthermore, $a_{ii}$ is supposed to be zero for all $i \in \{1, \ldots, n\}$ and $\mathcal{N}_i = \{v_j | \varepsilon_{ij} \in \varepsilon\}$ is the set of neighbors of agent $i$ [20]. $L_{n \times n}$ is Laplacian matrix of graph $\mathcal{G}$ is as follows

$$L = (l_{ij}), \quad l_{ij} = \begin{cases} \sum_{p \in \mathcal{N}_i} a_{ip}, & i = j \\ -a_{ij}, & i \neq j \end{cases}. \quad (1)$$

### C. Homotopy Method

Homotopy is the main idea for solving bilinear matrix inequalities (BMIs) to design the proposed controller in section 4. Assume solving $F(\zeta) = 0$ where $F(\zeta): \mathcal{R}^n \to \mathcal{R}^n$ is a nonlinear function. If $\hat{\zeta}$ be a solution of the mentioned problem, the main challenge of finding $\hat{\zeta}$ by iterative methods is to choose a starting point $\zeta_0$ such that it lies on the region of attraction for $\hat{\zeta}$ which guarantees the convergence of the solving process. Homotopy method eliminates the difficulty of finding an appropriate $\zeta_0$ for high dimensional problems.

In this method we consider a Homotopy function $H(\zeta, \eta): \mathcal{R}^n \times [0,1] \to \mathcal{R}^n$, where $H(\zeta, 1) = F(\zeta)$ and $H(\zeta, 0)$ is a function to calculate the starting point, such that $\zeta_0$ can be easily obtained by solving $H(\zeta_0, 0) = 0$. The desired solution $\hat{\zeta}$ is eventually achieved by increasing $\eta$ from 0 to 1 with sufficiently small steps; i.e. the solution of $H(\zeta, 1) = 0$ is equivalent to solving $F(\zeta) = 0$.

## III. MAIN RESULT

Assume an uncertain FOMAS consisting of $N$ agents, where the i-th agent dynamics are as follows

$$\begin{aligned} D^\alpha x_i(t) &= (\tilde{A} + \Delta \tilde{A}_i) x_i(t) + \tilde{B}_i u_i(t), \quad 0 < \alpha < 2 \\ y_i(t) &= \tilde{C} x_i(t), \qquad i = 1, \ldots, N \end{aligned}, \quad (2)$$

where $x_i \in \mathcal{R}^n$ denotes the state vector, $u_i \in \mathcal{R}^l$ is the control input, $y_i \in \mathcal{R}^q$ is the output vector of agent $i$. Furthermore $\tilde{A} \in \mathcal{R}^{n \times n}$, $\tilde{B}_i \in \mathcal{R}^{n \times l}$, and $\tilde{C}_i \in \mathcal{R}^{q \times n}$ are known constant matrices with appropriate dimensions. $\Delta \tilde{A}_i$ is time-invariant matrix with parametric uncertainty of the following form (see for example [21], [22])

$$\Delta \tilde{A}_i = \tilde{R} \delta_i(\sigma) \tilde{N}, \quad (3)$$
$$\delta_i(\sigma) = Z_i(\sigma)[I + JZ_i(\sigma)]^{-1}, \quad (4)$$
$$Sym\{J\} > 0, \quad (5)$$

where $\tilde{R} \in \mathcal{R}^{n \times m_0}$, $\tilde{N} \in \mathcal{R}^{m_0 \times n}$, and $J \in \mathcal{R}^{m_0 \times m_0}$ are real known matrices. The uncertain matrix $Z_i(\sigma) \in \mathcal{R}^{m_0 \times m_0}$ satisfies

$$Sym\{Z_i(\sigma)\} \geq 0, \quad (6)$$

where $\sigma \in \Omega$ with $\Omega$ being a compact set.

**Remark 1.** Condition (5) guarantees that $I + JZ_i(\sigma)$ is invertible for all $Z_i(\sigma)$ satisfying (6). Therefore $\delta_i(\sigma)$ in (4) is well defined ([23]).

In order to study the stability of fractional multi-agent systems and obtain main results we need following lemmas.

**Lemma 1.** ([24]) Let $A \in \mathcal{R}^{n \times n}$, $0 < \alpha < 1$ and $\theta = \alpha \pi / 2$. The fractional order system $D^\alpha x(t) = Ax(t)$ is asymptotically stable if and only if there exist two real symmetric positive definite matrices $P_{k1} \in \mathcal{R}^{n \times n}$, $k = 1, 2$, and two skew-symmetric matrices $P_{k2} \in \mathcal{R}^{n \times n}$, $k = 1, 2$, such that

$$\sum_{i=1}^{2} \sum_{j=1}^{2} Sym\{\Theta_{ij} \otimes (AP_{ij})\} < 0$$

$$\begin{bmatrix} P_{11} & P_{12} \\ -P_{12} & P_{11} \end{bmatrix} > 0 \quad \begin{bmatrix} P_{21} & P_{22} \\ -P_{22} & P_{21} \end{bmatrix} > 0, \quad (7)$$

where

$$\Theta_{11} = \begin{bmatrix} \sin\theta & -\cos\theta \\ \cos\theta & \sin\theta \end{bmatrix}, \quad \Theta_{12} = \begin{bmatrix} \cos\theta & \sin\theta \\ -\sin\theta & \cos\theta \end{bmatrix}$$

$$\Theta_{21} = \begin{bmatrix} \sin\theta & \cos\theta \\ -\cos\theta & \sin\theta \end{bmatrix}, \quad \Theta_{22} = \begin{bmatrix} -\cos\theta & \sin\theta \\ -\sin\theta & -\cos\theta \end{bmatrix} \quad (8)$$

.

The multi-agent system (2) can be represented in the following augmented form

$$\begin{aligned} D^\alpha x(t) &= \bar{A}_{\Delta,N} x(t) + Bu(t), \\ y(t) &= Cx(t), \end{aligned} \quad (9)$$

where $x(t) = [x_1(t)^T, \ldots, x_N(t)^T]^T$, $u(t) = [u_1(t)^T, \ldots, u_N(t)^T]^T$, and $y(t) = [y_1(t)^T, \ldots, y_N(t)^T]^T$ are general pseudo state, input, and output vectors respectively, and also $\bar{A}_{\Delta,N} = A_N + \overline{\Delta A}_N$ where $A_N = I_N \otimes \tilde{A}$ and $\overline{\Delta A}_N = diag(\Delta \tilde{A}_i)$, $B = diag(\tilde{B}_1, \tilde{B}_2, \ldots, \tilde{B}_N)$, and $C_N = I_N \otimes \tilde{C}$ are constant known matrices. Our aim is to define a fractional-order decentralized dynamic output feedback controller that guarantees consensus of agents.

**Definition 1.** System (2) achieves consensus asymptotically if the following condition holds for any $x_i(0) = x_{i0}$

$$\lim_{t \to \infty} |x_i(t) - x_j(t)| = 0, \quad i, j = 1, \ldots, N, i \neq j. \quad (10)$$

In order to solve the consensus problem for system (2), we use the following control protocol

$$D^\alpha x_{ci}(t) = A_{ci} x_{ci}(t) + B_{ci} \big( l_{ii} y_i(t) + \sum_{p \in \mathcal{N}_i} l_{ip} y_p(t) \big), \quad (11)$$

$u_i = C_{ci}x_i + D_{ci}\left(l_{ii}y_i(t) + \sum_{p \in \mathcal{N}_i} l_{ip}y_p(t)\right)$ ,
$i = 1, \ldots, N$.
Defining $x_c = [x_{c1}^T \ldots x_{cN}^T]^T$ as controller pseudo state vector, where $x_{c1}, \ldots, x_{cN} \in \mathcal{R}^{n_c}$ in which $n_c$ is the controller order, $A_c = diag(A_{c1}, A_{c2}, \ldots, A_{cN})$, $B_c = diag(B_{c1}, \ldots, B_{cN})$, $C_c = diag(C_{c1}, \ldots, C_{cN})$, $D_c = diag(D_{c1}, \ldots, D_{cN})$ as controller matrices and $\mathcal{N}_i$ denotes the neighbors of agent $i$. To simplify (11) we can represent the control protocol in the following form

$$D^\alpha x_c(t) = A_c x_c(t) + B_c L_q y(t), \qquad (12)$$
$$u = C_c x_c + D_c L_q y(t),$$

where $L_q = L \otimes I_q$, in which $L$ is the Laplacian matrix of corresponding graph of FOMAS. Implementing control protocol (12) to the main system (9) yields

$$D^\alpha \begin{bmatrix} x(t) \\ x_c(t) \end{bmatrix} = \begin{bmatrix} \bar{A}_{\Delta,N} + BD_c L_q C_N & BC_c \\ B_c L_q C_N & A_c \end{bmatrix} \begin{bmatrix} x(t) \\ x_c(t) \end{bmatrix} = \left( \begin{bmatrix} A_N + BD_c L_q C_N & BC_c \\ B_c L_q C_N & A_c \end{bmatrix} + \begin{bmatrix} \overline{\Delta A}_N & 0 \\ 0 & 0 \end{bmatrix} \right) \begin{bmatrix} x(t) \\ x_c(t) \end{bmatrix}. \qquad (13)$$

**Lemma 4.** [11] $L_q C_N = C_N L_n$.
Proof. $L_q C_N = (L \otimes I_q)(I_N \otimes \tilde{C}) = (LI_N) \otimes (I_q \tilde{C}) = (I_N L) \otimes (\tilde{C} I_n) = (I_N \otimes \tilde{C})(L \otimes I_n) = C_N L_n$.
**Lemma 5.** [11] $L_n A_N = A_N L_n$.
Proof. $L_n A_N = (L \otimes I_n)(I_N \otimes \tilde{A}) = (LI_N) \otimes (I_n \tilde{A}) = (I_N L) \otimes (\tilde{A} I_n) = (I_N \otimes \tilde{A})(L \otimes I_n) = A_N L_n$.

**Theorem 1.** The output feedback controller (11) solves the consensus problem of the system (2) with $0 < \alpha < 1$, if there exist two real symmetric positive definite matrices $X_{k1} \in \mathcal{R}^{n \times n}$, $k = 1,2$, and two skew-symmetric matrices $X_{k2} \in \mathcal{R}^{n \times n}$, $k = 1,2$, and a real constant $\mu > 0$ such that

$$\sum_{i=1}^{2}\sum_{j=1}^{2} \begin{bmatrix} \Xi_{X,ij} & \Theta_{ij} \otimes M & I_2 \otimes (X_{ij}^T N^T) \\ \star & -\mu I & \mu I \\ \star & \star & \Upsilon \end{bmatrix} < 0, \qquad (14)$$

where
$\Xi_{X_{ij}} = Sym\{\Theta_{ij} \otimes (A_{cl} X_{ij})\}$,
$\Upsilon = -Sym\{I_2 \otimes (I_{N-1} \otimes J)\} - \mu I$,
$A_{cl} = \begin{bmatrix} A_{N-1} + (\hat{L} \otimes I_n)BD_c C_N L_n \hat{L}_n^\uparrow & (\hat{L} \otimes I_n)BC_c \\ B_c C_N L_n \hat{L}_n^\uparrow & A_c \end{bmatrix}$, (15)
$M = [I_{N-1} \otimes \widetilde{M}^T \mid \mathbf{0}_{(N-1).m_0 \times N.n_c}]^T$,
$R = [I_{N-1} \otimes \tilde{R} \mid \mathbf{0}_{(N-1).m_0 \times N.n_c}]$.
in which $\hat{L} \in \mathcal{R}^{(N-1) \times N}$ is obtained by removing a dependent row of $L$. And
$\Theta_{11} = \begin{bmatrix} \sin\theta & -\cos\theta \\ \cos\theta & \sin\theta \end{bmatrix}$, $\Theta_{12} = \begin{bmatrix} \cos\theta & \sin\theta \\ -\sin\theta & \cos\theta \end{bmatrix}$,
$\Theta_{21} = \begin{bmatrix} \sin\theta & \cos\theta \\ -\cos\theta & \sin\theta \end{bmatrix}$, $\Theta_{22} = \begin{bmatrix} -\cos\theta & \sin\theta \\ -\sin\theta & -\cos\theta \end{bmatrix}$,
$\theta = \alpha\pi/2$.

**Proof.** The idea is to convert the consensus problem of (13) into a stabilization one. In order to guarantee the consensus of FOMAS (13) via stabilization problem, we define a new system with the following states

$$\bar{x} = L_n x, \qquad L_n = L \otimes I_n. \qquad (16)$$

It is obvious that $\bar{x}$ provides a relation between $x_i$ and its neighbors. The stability of the transformed system ensures the convergence of relative states difference of system (13) to zero which is equivalent to the consensus of system (13) defined in (10). Nevertheless, the transformed system by (16) has a redundancy due to rank deficiency of $L_n$. Removing a row of the matrix $L$ eliminates the system redundancy, so we modify (16) as follows

$$x_r = \hat{L}_n x, \qquad \hat{L}_n = \hat{L} \otimes I_n. \qquad (17)$$

In order to express the pseudo state space representation of system with respect to $x_r$ we have

$$D^\alpha {}_N\Lambda^n_{n_c} \begin{bmatrix} x(t) \\ x_c(t) \end{bmatrix} = {}_N\Lambda^n_{n_c} \begin{bmatrix} \bar{A}_{\Delta,N} + BD_c L_q C_N & BC_c \\ B_c L_q C_N & A_c \end{bmatrix} {}_N\Lambda^{n\,\uparrow}_{n_c} {}_N\Lambda^n_{n_c} \begin{bmatrix} x(t) \\ x_c(t) \end{bmatrix}, \qquad (18)$$

where
$${}_N\Lambda^n_{n_c} = diag(\hat{L}_n, I_{N.n_c}) \qquad (19)$$

considering
$A_{\Delta,N} = I_N \otimes \Delta\tilde{A}$,
$\Delta\tilde{A} = \tilde{R}\delta(\sigma)\tilde{N}$, (20)
$\delta = \{\delta_j \mid \|\delta_j\| > \|\delta_i\|, i \neq j\}$, $i = 1, \ldots, N$,
some calculations yield the closed loop system

$$D^\alpha \begin{bmatrix} x_r(t) \\ x_c(t) \end{bmatrix} = \begin{bmatrix} \hat{L}_n A_{\Delta,N} \hat{L}_n^\uparrow + \hat{L}_n BD_c L_q C_N \hat{L}_n^\uparrow & \hat{L}_n BC_c \\ B_c L_q C_N \hat{L}_n^\uparrow & A_c \end{bmatrix} \begin{bmatrix} x_r(t) \\ x_c(t) \end{bmatrix}. \qquad (21)$$

Lemma 4 and Lemma 5 help us to match the matrices dimensions in the resulting system matrix. Using these lemmas, it can be easily obtained that $\hat{L}_n A_{\Delta,N} \hat{L}_n^\uparrow = I_{N-1} \otimes (\tilde{A} + \Delta\tilde{A}) = A_{\Delta,N-1}$ and $L_q C_N \hat{L}_n^\uparrow = C_N L_n \hat{L}_n^\uparrow$. Then, defining $C_r = C_N L_n \hat{L}_n^\uparrow$, the closed loop system is achieved as follows

$$D^\alpha \begin{bmatrix} x_r(t) \\ x_c(t) \end{bmatrix} = \begin{bmatrix} A_{\Delta,N-1} + \hat{L}_n BD_c C_r & \hat{L}_n BC_c \\ B_c C_r & A_c \end{bmatrix} \begin{bmatrix} x_r(t) \\ x_c(t) \end{bmatrix}, \qquad (22)$$
$A_{cl,\Delta} = A_{cl} + \Delta A_{cl}$,
where
$A_{cl} = \begin{bmatrix} A_{N-1} + \hat{L}_n BD_c C_r & \hat{L}_n BC_c \\ B_c C_r & A_c \end{bmatrix}$, (23)
$\Delta A_{cl} = M[I_{(N-1).n} \otimes \delta(\sigma)]R$,
where
$M = [I_{N-1} \otimes \widetilde{M}^T \mid \mathbf{0}_{(N-1).m_0 \times N.n_c}]^T$, (24)
$R = [I_{N-1} \otimes \tilde{R} \mid \mathbf{0}_{(N-1).m_0 \times N.n_c}]$.

To summarize the above discussion we have reduced the consensus problem to the stabilization of the closed loop system (22). It follows from Lemma 1 that the system (22) is stable if and only if there exist two real symmetric positive definite matrices $P_{k1} \in \mathcal{R}^{n \times n}$, $k = 1,2$, and two skew-symmetric matrices $P_{k2} \in \mathcal{R}^{n \times n}$, $k = 1,2$, such that
$\sum_{i=1}^{2}\sum_{j=1}^{2} Sym\{\Theta_{ij} \otimes (A_{cl,\Delta} P_{ij})\} < 0$,
$\begin{bmatrix} P_{11} & P_{12} \\ -P_{12} & P_{11} \end{bmatrix} > 0$, $\begin{bmatrix} P_{21} & P_{22} \\ -P_{22} & P_{21} \end{bmatrix} > 0$,
where
$\Theta_{11} = \begin{bmatrix} \sin\theta & -\cos\theta \\ \cos\theta & \sin\theta \end{bmatrix}$, $\Theta_{12} = \begin{bmatrix} \cos\theta & \sin\theta \\ -\sin\theta & \cos\theta \end{bmatrix}$,
$\Theta_{21} = \begin{bmatrix} \sin\theta & \cos\theta \\ -\cos\theta & \sin\theta \end{bmatrix}$, $\Theta_{22} = \begin{bmatrix} -\cos\theta & \sin\theta \\ -\sin\theta & -\cos\theta \end{bmatrix}$,
$\theta = \alpha\pi/2$.

$\sum_{i=1}^{2}\sum_{j=1}^{2} Sym\{\Theta_{ij} \otimes (A_{cl,\Delta} P_{ij})\} < 0 \Leftrightarrow$
$\sum_{i=1}^{2}\sum_{j=1}^{2} Sym\{\Theta_{ij} \otimes ((A_{cl} + \Delta A_{cl})P_{ij})\} < 0 \Leftrightarrow$
$\sum_{i=1}^{2}\sum_{j=1}^{2} Sym\{\Theta_{ij} \otimes (A_{cl}P_{ij})\} +$
$\sum_{i=1}^{2}\sum_{j=1}^{2} Sym\{\Theta_{ij} \otimes (\Delta A_{cl} P_{ij})\} < 0 \Leftrightarrow$ (25)
$\sum_{i=1}^{2}\sum_{j=1}^{2} Sym\{\Theta_{ij} \otimes (A_{cl}P_{ij})\} +$
$\sum_{i=1}^{2}\sum_{j=1}^{2} Sym\{\Theta_{ij} \otimes (M(I_{N-1} \otimes \delta(\sigma))RP_{ij})\} <$

$0 \Leftrightarrow \sum_{i=1}^{2}\sum_{j=1}^{2} Sym\{\Theta_{ij} \otimes (A_{cl}P_{ij})\} +$
$\sum_{i=1}^{2}\sum_{j=1}^{2} Sym\{\Theta_{ij} \otimes (M\Delta(\sigma)RP_{ij})\} < 0 \Leftrightarrow \Psi +$
$\sum_{i=1}^{2}\sum_{j=1}^{2} Sym\{\widehat{M}_{ij}\widehat{\Delta}(\sigma)\widehat{R}\widehat{P}_{ij}\} < 0,$
where $\Psi = \sum_{i=1}^{2}\sum_{j=1}^{2} Sym\{\Theta_{ij} \otimes (A_{cl}P_{ij})\}$, $\widehat{M}_{ij} = \Theta_{ij} \otimes M$, $\widehat{R} = I_2 \otimes R$, $\widehat{Z}(\sigma) = I_2 \otimes (I_{N-1} \otimes Z(\sigma))$, $\hat{J} = I_2 \otimes (I_{N-1} \otimes J)$, $\widehat{\Delta}(\sigma) = I_2 \otimes (I_{N-1} \otimes \delta(\sigma)) = \widehat{Z}(\sigma)[I - \hat{J}\widehat{Z}(\sigma)]^{-1}$, $\widehat{P}_{ij} = I_2 \otimes P_{ij} (i,j = 1,2)$.

Define
$\widehat{W} = Sym(\hat{J})$, $\hat{S} = \sum_{i=1}^{2}\sum_{j=1}^{2} \widehat{P}_{ij}$,
$\widehat{Q}_{ij} = \widehat{W}^{-1/2}(\varepsilon^{-1}\widehat{M}_{ij}^T + \varepsilon\widehat{R}\widehat{P}_{ij}) -$  (26)
$\widehat{W}^{1/2}\widehat{\Delta}^T(\sigma)\varepsilon^{-1}\widehat{M}_{ij}^T.$

It follows from Lemma 3 that $Sym\{\widehat{\Delta}(\sigma)\} - \widehat{\Delta}(\sigma)\widehat{W}\widehat{\Delta}^T(\sigma) > 0$, and for any $\varepsilon > 0$ the following inequality holds
$-\sum_{i=1}^{2}\sum_{j=1}^{2} \widehat{Q}_{ij}^T\widehat{Q}_{ij} \leq 0 \Leftrightarrow$
$\sum_{i=1}^{2}\sum_{j=1}^{2} \left[ -Sym\{\widehat{M}_{ij}\widehat{W}^{-1}\widehat{R}\widehat{P}_{ij}\} - \varepsilon^{-2}\widehat{M}_{ij}\widehat{W}^{-1}\widehat{M}_{ij}^T - \varepsilon^2\widehat{P}_{ij}^T\widehat{R}^T\widehat{W}^{-1}\widehat{R}\widehat{P}_{ij} + Sym\{\widehat{M}_{ij}\widehat{\Delta}(\sigma)\widehat{R}\widehat{P}_{ij}\} + \varepsilon^{-2}\widehat{M}_{ij}\left(Sym\{\widehat{\Delta}(\sigma)\} - \widehat{\Delta}(\sigma)\widehat{W}\widehat{\Delta}^T(\sigma)\right)\widehat{M}_{ij}^T \right] \leq 0 \Rightarrow$  (27)
$\sum_{i=1}^{2}\sum_{j=1}^{2} [-Sym\{\widehat{M}_{ij}\widehat{W}^{-1}\widehat{R}\widehat{P}_{ij}\} - \varepsilon^{-2}\widehat{M}_{ij}\widehat{W}^{-1}\widehat{M}_{ij}^T - \varepsilon^2\widehat{P}_{ij}^T\widehat{R}^T\widehat{W}^{-1}\widehat{R}\widehat{P}_{ij}^T] +$
$\sum_{i=1}^{2}\sum_{j=1}^{2} Sym\{\widehat{M}_{ij}\widehat{\Delta}(\sigma)\widehat{R}\widehat{P}_{ij}\} \leq 0.$

It follows from inequality (27) that inequality (25) holds if
$\Psi + \sum_{i=1}^{2}\sum_{j=1}^{2} [-Sym\{\widehat{M}_{ij}\widehat{W}^{-1}\widehat{R}\widehat{P}_{ij}\} - \varepsilon^{-2}\widehat{M}_{ij}\widehat{W}^{-1}\widehat{M}_{ij}^T - \varepsilon^2\widehat{P}_{ij}^T\widehat{R}^T\widehat{W}^{-1}\widehat{R}\widehat{P}_{ij}^T] < 0,$  (28)

inequality (28) is equivalent to that there exist $\varepsilon > 0$ and $\mu > 0$ such that
$\Psi + \sum_{i=1}^{2}\sum_{j=1}^{2} [-Sym\{\widehat{M}_{ij}\widehat{W}^{-1}\widehat{R}\widehat{P}_{ij}\} - \varepsilon^{-2}\widehat{M}_{ij}(\widehat{W}^{-1} + \mu^{-1}I)\widehat{M}_{ij}^T - \varepsilon^2\widehat{P}_{ij}^T\widehat{R}^T\widehat{W}^{-1}\widehat{R}\widehat{P}_{ij}^T] < 0,$  (29)

which is equivalent to that there exist $\varepsilon > 0$ and $\mu > 0$ such that
$\Psi +$
$\sum_{i=1}^{2}\sum_{j=1}^{2} \begin{bmatrix} \varepsilon^{-1}\widehat{M}_{ij}^T \\ \varepsilon\widehat{R}\widehat{P}_{ij} \end{bmatrix}^T \begin{bmatrix} \widehat{W}^{-1} + \mu^{-1}I & \widehat{W}^{-1} \\ \widehat{W}^{-1} & \widehat{W}^{-1} \end{bmatrix} \begin{bmatrix} \varepsilon^{-1}\widehat{M}_{ij}^T \\ \varepsilon\widehat{R}\widehat{P}_{ij} \end{bmatrix}$
$= \Psi +$  (30)
$\sum_{i=1}^{2}\sum_{j=1}^{2} \begin{bmatrix} \varepsilon^{-1}\widehat{M}_{ij}^T \\ \varepsilon\widehat{R}\widehat{P}_{ij} \end{bmatrix}^T \begin{bmatrix} \mu I & -\mu I \\ -\mu I & \widehat{W} + \mu I \end{bmatrix}^{-1} \begin{bmatrix} \varepsilon^{-1}\widehat{M}_{ij}^T \\ \varepsilon\widehat{R}\widehat{P}_{ij} \end{bmatrix} < 0.$

Applying Schur complement, inequality (30) become equivalent to
$\sum_{i=1}^{2}\sum_{j=1}^{2} \begin{bmatrix} \Xi_{P_{ij}} & \varepsilon^{-1}\widehat{M}_{ij} & \varepsilon\widehat{P}_{ij}^T\widehat{R}^T \\ \star & -\mu I & \mu I \\ \star & \star & -(\widehat{W} + \mu I) \end{bmatrix} < 0.$  (31)

Where
$\Xi_{P,ij} = Sym\{\Theta_{ij} \otimes (A_{cl}P_{ij})\}$  (32)

Pre- and post-multiplying the left side of the inequality (31) by $diag(\varepsilon I, I, I)$, results
$\sum_{i=1}^{2}\sum_{j=1}^{2} \begin{bmatrix} \varepsilon^2\Xi_{P_{ij}} & \widehat{M}_{ij} & \varepsilon\widehat{P}_{ij}^T\widehat{R}^T \\ \star & -\mu I & \mu I \\ \star & \star & -(\widehat{W} + \mu I) \end{bmatrix} < 0, \Leftrightarrow$  (33)

$\sum_{i=1}^{2}\sum_{j=1}^{2} \begin{bmatrix} \varepsilon^2\Xi_{X_{ij}} & \Theta_{ij} \otimes M & I_2 \otimes (\varepsilon^2 P_{ij}^T R^T) \\ \star & -\mu I & \mu I \\ \star & \star & \Upsilon \end{bmatrix} < 0.$

Defining $X_{ij} = \varepsilon^2 P_{ij}, i,j = 1,2$, the inequality (33) becomes equivalent to (14). This completes the proof. ∎

**Corollary 1.** Consider the multi-agent system (2) with $0 < \alpha < 1$, without the uncertainty. the output feedback controller (11) solves the consensus problem, if there exist two real symmetric positive definite matrices $X_{k1} \in \mathcal{R}^{n \times n}$, $k = 1,2$, and two skew-symmetric matrices $X_{k2} \in \mathcal{R}^{n \times n}$, $k = 1,2$, such that

$\sum_{i=1}^{2}\sum_{j=1}^{2} Sym\{\Theta_{ij} \otimes (A_{cl}X_{ij})\} < 0,$
$\begin{bmatrix} X_{11} & X_{12} \\ -X_{12} & X_{11} \end{bmatrix} > 0, \quad \begin{bmatrix} X_{21} & X_{22} \\ -X_{22} & X_{21} \end{bmatrix} > 0,$  (34)

where
$A_{cl} = \begin{bmatrix} I_{N-1} \otimes \tilde{A} + (\hat{L} \otimes I_n)BD_cC_NL_n\hat{L}_n^{\uparrow} & (\hat{L} \otimes I_n)BC_c \\ B_cC_NL_n\hat{L}_n^{\uparrow} & A_c \end{bmatrix},$  (35)

in which $\hat{L} \in \mathcal{R}^{(N-1) \times N}$ is obtained by removing a dependent row of $L$. And
$\Theta_{11} = \begin{bmatrix} \sin\theta & -\cos\theta \\ \cos\theta & \sin\theta \end{bmatrix}, \quad \Theta_{12} = \begin{bmatrix} \cos\theta & \sin\theta \\ -\sin\theta & \cos\theta \end{bmatrix},$
$\Theta_{21} = \begin{bmatrix} \sin\theta & \cos\theta \\ -\cos\theta & \sin\theta \end{bmatrix}, \quad \Theta_{22} = \begin{bmatrix} -\cos\theta & \sin\theta \\ -\sin\theta & -\cos\theta \end{bmatrix},$

Proof. Set $\Delta\tilde{A} = 0$ in (22). The proof procedure is similar to that of Theorem 1, and it is omitted. ∎

**Remark 2.** In Theorem 1 and Corollary 1 the obtained dynamic output feedback controllers can be reduced to static ones by putting $n_c = 0$.

## IV. DISCUSSION ON SOLVING INEQUALITY (14)

Theorem 1 gives sufficient condition for solving the consensus problem of FOMAS (9) for $0 < \alpha < 1$. The inequality (14) is bilinear matrix inequality (BMI), since matrix $A_{cl}$ containing varying terms, is multiplied by $X_{ij}$, ($i,j = 1,2$). The closed loop matrix (22) is rewritten in the following specified form in order to simplify the solving process of the inequality (14).

$A_{cl,\Delta} = A_K + B_K K C_K,$  (36)
where
$A_K = \begin{bmatrix} A_{\Delta,N-1} & 0 \\ 0 & 0 \end{bmatrix}, \quad B_K = \begin{bmatrix} \hat{L}_n B & 0 \\ 0 & I \end{bmatrix}, \quad K = \begin{bmatrix} D_c & C_c \\ B_c & A_c \end{bmatrix}, \quad C_K = \begin{bmatrix} C_r & 0 \\ 0 & I \end{bmatrix}.$  (37)

Since the inequality (14) can't be converted into LMI form, and consequently it is difficult to solve, we use the Homotopy based approach proposed in [25], to deal with this obstacle.

By introducing $\eta$ as a real number varying between 0 and 1, the matrix function $F$ is defined as follows
$F(X_{ij}, K, \eta) = (1-\eta)F_1(X_{ij}) + \eta F_2(X_{ij}, K), \quad i,j = 1,2.$  (38)

Where $F_1(X_{ij})$ and $F_2(X_{ij}, K)$ are presented in (39) with $\sigma \geq 0$, in which $Q$ is a matrix with suitable dimensions such that $(A_K - Q)$ is Hurwitz [25]. Since

$$F(X_{ij}, K, \eta) = \begin{cases} F_1(X_{ij}), & \eta = 0 \\ F_2(X_{ij}, K), & \eta = 1 \end{cases}, \quad (40)$$

finding a feasible solution to inequality (14) is equivalent to solving the following parametrized problems

$$F(X_{ij}, K, \eta) < 0, \quad \eta \in [0,1]. \quad (41)$$

by gradually increasing $\eta$ from 0 to 1. For $\eta = 0$, the inequality (41) becomes a Linear Matrix Inequality (LMI), which can be solved by appropriate Solvers and Parsers. Then we solve the inequality (41) considering $X_{ij}$ and $K$ to be fixed alternatively where $\eta$ has been increased during this procedure by supposing $\eta = i/Z$ with $i = 1, 2, ..., T$ and a large positive integer $T$.

## V. SIMULATIONS

In this section, the proposed example demonstrates the effectiveness of the designed decentralized dynamic output feedback controller for the consensus of a fractional-order multi-agent system. Various solvers and parsers can be utilized to determine variables satisfying feasibility problem. In this paper simulation results are obtained using YALMIP parser [26], implemented as a toolbox in Matlab [27].

### A. Numerical example

A connected network with four agents is considered as shown in Figure 1. The dynamics of each agent are represented as follows

$$D^\alpha x_i(t) = (\tilde{A} + \Delta\tilde{A}_i) x_i(t) + B_i U_{di}(t),$$
$$y_i(t) = \tilde{C} x_i(t), \quad i = 1, 2, 3, \quad (42)$$

where

$$\tilde{A} = \begin{bmatrix} -1 & 1 & 0 & 0 \\ 1 & -3 & 0 & 1 \\ 0 & 0 & 0 & 1 \\ 0 & 0 & -1 & 0 \end{bmatrix}, B_i = \begin{bmatrix} 1 \\ 1 \\ 0 \\ 1 \end{bmatrix}, \tilde{C} = \begin{bmatrix} 1 \\ 0 \\ 1 \\ 0 \end{bmatrix}^T,$$
$$\tilde{M} = [0.2 \ 0 \ -0.1 \ 0.3]^T, \tilde{R} = [0 \ 0.2 \ 0.4 \ -0.2], J = 1,$$
$$\delta_1 = 0.5, \delta_2 = -0.4, \delta_3 = 0.1, \delta_4 = 0.8, \quad (43)$$

where $\delta = 0.8$. The corresponding Laplacian matrix of the mentioned network is $L = \begin{bmatrix} 2 & -1 & 0 & -1 \\ -1 & 2 & -1 & 0 \\ 0 & -1 & 2 & -1 \\ -1 & 0 & -1 & 2 \end{bmatrix}$.

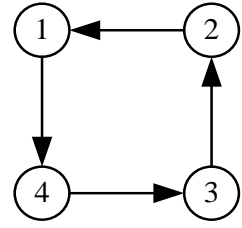

Figure 1. Network topology o FOMAS.

This example is solved for $\alpha = 0.8$, where the resulted controllers obtained by Theorem 1 are presented in Table 1. For more efficient convergence of solving (38), $\eta = 0.1$ and $T = 10$ are chosen.

Figure 2 illustrates the states trajectories of agents for initial states $x_0 = [x_{10} \ x_{20} \ x_{30} \ x_{40}]$ using the output feedback controllers presented in Table 1, where

$$x_{10} = \begin{bmatrix} 2.5 \\ -3 \\ -2.5 \\ 1.5 \end{bmatrix}, x_{20} = \begin{bmatrix} 0.1 \\ 0.5 \\ 2 \\ 1 \end{bmatrix}, x_{30} = \begin{bmatrix} 0.1 \\ 0.8 \\ 0.54 \\ 1.1 \end{bmatrix}, x_{40} = \begin{bmatrix} -3 \\ -2.5 \\ -2.3 \\ 3 \end{bmatrix}, \quad (44)$$

are the initial states of each agent. Obviously, the corresponding trajectories of all agents asymptotically reach agreement.

The consensus error of proposed static and dynamic controllers' indices are summarized in Table 2 where the results indicate that increasing controller order reduces the settling time of the controller effort. Since the controller objective is the consensus of the corresponding states of agents, the vanishing of the consensus error, indicates that the consensus of agents is achieved. Thus, there is a direct relationship between the controller order and the rate of the consensus. Besides increasing the controller order leads to a

$$F_1(X_{ij}) = \sum_{i=1}^{2}\sum_{j=1}^{2} \begin{bmatrix} Sym\{\Theta_{ij} \otimes ((A_K - Q)X_{ij})\} & \Theta_{ij} \otimes M & I_2 \otimes (X_{ij}^T R^T) \\ \star & -\mu I & \mu I \\ \star & \star & \Upsilon \end{bmatrix},$$

$$F_2(X_{ij}, K) = \sum_{i=1}^{2}\sum_{j=1}^{2} \begin{bmatrix} \Xi_{X,ij} & \Theta_{ij} \otimes M & I_2 \otimes (X_{ij}^T R^T) \\ \star & -\mu I & \mu I \\ \star & \star & \Upsilon \end{bmatrix}, \quad (39)$$

Table 1. Controller parameters obtained by (I) Theorem 1 and (II) Corollary 1.

| Order of controller | | | Agent 1 | Agent 2 | Agent 3 | Agent 4 |
|---|---|---|---|---|---|---|
| **0** | $D_c$ | (I) | $-2.2501$ | $-3.0239$ | $-2.2501$ | $-1.3580$ |
| | | (II) | $-2.7534$ | $-3.1985$ | $-2.7534$ | $-1.8184$ |
| **2** | $A_c$ | (I) | $\begin{bmatrix} -53.7197 & -15.1371 \\ 22.4322 & -24.3364 \end{bmatrix}$ | $\begin{bmatrix} -35.3656 & 39.5514 \\ 13.0965 & -25.5037 \end{bmatrix}$ | $\begin{bmatrix} -43.6527 & 8.7360 \\ 11.9074 & -25.4178 \end{bmatrix}$ | $\begin{bmatrix} -36.5995 & 10.9339 \\ 12.8142 & -29.2978 \end{bmatrix}$ |
| | | (II) | $\begin{bmatrix} -36.1669 & 6.7979 \\ 6.7979 & -36.1669 \end{bmatrix}$ | $\begin{bmatrix} -37.0864 & 5.8784 \\ 5.8784 & -37.0864 \end{bmatrix}$ | $\begin{bmatrix} -36.1669 & 6.7979 \\ 6.7979 & -36.1669 \end{bmatrix}$ | $\begin{bmatrix} -36.9240 & 6.0408 \\ 6.0408 & -36.9240 \end{bmatrix}$ |
| | $B_c$ | (I) | $[43.6203 \ -1.8477]^T$ | $[-10.2150 \ -4.7566]^T$ | $[16.6661 \ 4.5128]^T$ | $[7.3644 \ 5.2846]^T$ |
| | | (II) | $[0.8343 \ 0.8343]^T$ | $[1.2489 \ 1.2489]^T$ | $[0.8343 \ 0.8343]^T$ | $[-0.8981 \ -0.8981]^T$ |
| | $C_c$ | (I) | $[-2.1540 \ -4.7027]$ | $[6.5429 \ -7.9858]$ | $[-6.0084 \ -4.1185]$ | $[3.3085 \ 3.0719]$ |
| | | (II) | $[0.2524 \ 0.2524]$ | $[0.3631 \ 0.3631]$ | $[0.2524 \ 0.2524]$ | $[-1.3966 \ -1.3966]$ |
| | $D_c$ | (I) | $-10.0066$ | $-10.0034$ | $-9.9941$ | $-8.8604$ |
| | | (II) | $-3.4923$ | $-4.9571$ | $-3.4923$ | $-2.3834$ |

slight decrement in consensus error indices, and more quick and efficient consensus.

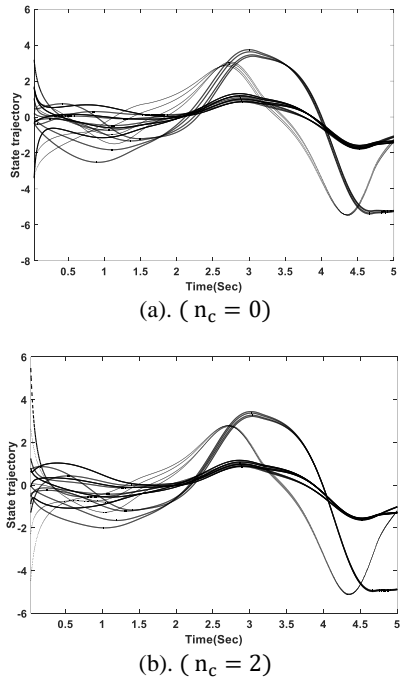

(a). ($n_c = 0$)

(b). ($n_c = 2$)

Figure 2. States trajectories of multi-agent system shown in Figure 1 using the proposed controllers in Theorem 1.

Table 2. Consensus error indices for method in Theorem 1.

| Controller order | | ISE | IAE | ITSE | ITAE |
|---|---|---|---|---|---|
| 0 | Agent 1 | 3.2706 | 1.0023 | 0.0984 | 0.4007 |
| | Agent 2 | 6.2087 | 2.7805 | 1.7648 | 3.5645 |
| | Agent 3 | 3.2706 | 0.5254 | 0.0293 | 0.2255 |
| | Agent 4 | 1.9227 | 1.3530 | 0.3807 | 1.7665 |
| 2 | Agent 1 | 0.5385 | 0.6690 | 0.1086 | 0.5835 |
| | Agent 2 | 0.1311 | 0.3972 | 0.0359 | 0.3099 |
| | Agent 3 | 1.2063 | 0.6203 | 0.0424 | 0.2753 |
| | Agent 4 | 0.5075 | 0.3052 | 0.0141 | 0.1723 |

## VI. CONCLUSION

In this paper decentralized static and dynamic output feedback controllers for the consensus of the fractional-order multi-agent systems are proposed. First, a new FOMAS with transformed states is defined, in which the stability of this new system is equivalent to the consensus of the main system. Second, sufficient conditions for the stability of new system using fractional-order systems stability theorems and Schur complement are obtained in the form of matrix inequalities. Third, the controller unknown parameters obtained by solving matrix inequalities using a Homotopy based method. Eventually, some numerical examples are presented to illustrate the effectiveness of the proposed dynamic output feedback controller design methods for FOMASs.